\documentclass{amsart}
\usepackage{amsmath,amssymb,amsthm,url}
\newtheorem{thr}{Theorem}[section]
\newtheorem{lem}[thr]{Lemma}

\theoremstyle{definition}

\newtheorem{prob}[thr]{Problem}

\theoremstyle{remark}

\numberwithin{equation}{section}

\def\M{{\mathcal M}}
\def\ep{{\varepsilon}}
\def\V{{\mathcal V}}

\begin{document}

\title{Column subset selection is NP-complete}

\author{Yaroslav Shitov}
\address{National Research University Higher School of Economics, 20 Myasnitskaya Ulitsa, Moscow 101000, Russia}
\email{yaroslav-shitov@yandex.ru}



\begin{abstract}
Let $M$ be a real $r\times c$ matrix and let $k$ be a positive integer. In the \textit{column subset selection} problem (CSSP), we need to minimize the quantity $\|M-SA\|$, where $A$ can be an arbitrary $k\times c$ matrix, and $S$ runs over all $r\times k$ submatrices of $M$. This problem and its applications in numerical linear algebra are being discussed for several decades, but its algorithmic complexity remained an open issue. We show that CSSP is NP-complete.
\end{abstract}

\maketitle

\section{Introduction}

Throughout this paper, we denote by $\|M\|$ the \textit{Frobenius norm} of a real matrix $M$, that is, the square root of the sum of squares of the elements of $M$. If $M$ has $r$ rows and $c$ columns, then we write $M\in\mathbb{R}^{r\times c}$. For any integer $k\leqslant c$, we denote
\begin{equation}\label{eqdefd}
\delta_k(M)=\min\limits_{S,A}\|M-SA\|,
\end{equation}
where $A$ and $S$ run, respectively, over $\mathbb{R}^{k\times c}$ and over all $r\times k$ submatrices of $M$. If we fix $S$ in the definition of $\delta_k(M)$, then the optimal value of $A$ is known to equal $S^+M$, where $X^+$ denotes the pseudoinverse of $X$.

The column subset selection problem (CSSP) is the task of detecting a submatrix $S$ realizing the minimum in~\eqref{eqdefd}. The quantity $\delta_k(M)$ can be seen as a value of how close a matrix can be to the linear space spanned by a tuple of $k$ of its columns, which makes CSSP important for problems of low-rank matrix approximation. In fact, CSSP is being studied for more than twenty years (see~\cite{CH}), and many new areas of its application have arisen since that time. Besides the low-rank approximation theory~(\cite{BDM, BMD}), the subset selection problem arises naturally in statistical data analysis~(\cite{KPS}), large data analysis~(\cite{PT}), optimal experiment design~(\cite{dHM}), statistical computing~(\cite{CH}), artificial intelligence~(\cite{AXMS}), networking~(\cite{JB}), rank-deficient least squares computation~(\cite{FK, IKP}), and other branches of modern applied mathematics.

The CSSP problem has been widely discussed in the applied mathematics community, but its algorithmic complexity status remained unknown. In 2009, Boutsidis, Mahoney, and Drineas~(\cite{BMD}) pointed out that the NP-hardness of CSSP is an open problem, and this problem has been discussed in numerous sources since then~(\cite{AXMS, CMI, KS, KPNS, KPS2}). A recent paper~\cite{Civ} presents a first step towards a solution: \c{C}ivril proves that CSSP is UG-hard, that is, he gives a conditional NP-hardness proof subject to the validity of the so-called \textit{Unique Games conjecture}. Since there is currently no clear evidence of whether this conjecture is true or false (see~\cite{Khot}), the paper by \c{C}ivril left the complexity status of CSSP wide open.

\section{The result}

The goal of our paper is to determine the complexity status of CSSP and to prove that this problem is NP-hard. Moreover, we prove that the CSSP problem restricted to rational matrices is NP-complete in the Turing model of computation. This problem can be stated formally in the following way.

\begin{prob}\label{prob11} (COLUMN SUBSET SELECTION.)

\noindent Given: A matrix $M\in\mathbb{Q}^{r\times c}$, an integer $k\in\{1,\ldots,c\}$, and a positive rational threshold $\tau$.

\noindent Question: Is $\delta_k(M)\leqslant\sqrt{\tau}$?
\end{prob}

Recall that the optimal matrix $A$ in~\eqref{eqdefd} can be expressed as $S^+M$, where $X^+$ denotes the pseudoinverse of $X$. In fact, many authors prefer writing $SS^+M$ instead of $S A$ in~\eqref{eqdefd}, and this gives an equivalent definition of what we denote by $\delta_k(M)$. In particular, we see that Problem~\ref{prob11} belongs to NP: If we correctly choose the optimal submatrix $S$, then we are able to compute $\delta_k(M)$ in polynomial time. Our goal is to prove the following result.

\begin{thr}
Problem~\ref{prob11} is NP-complete.
\end{thr}

As said above, Problem~\ref{prob11} belongs to NP, so we only need to construct a polynomial transformation to Problem~\ref{prob11} from some known NP-complete problem. To this end, we use the classical problem of \textit{graph three-coloring}, see~\cite{Karp} for details. We consider simple graphs $G=(V,E)$, where $V$ is a finite non-empty set of vertices, and $E$ is a finite non-empty set of edges. A function $\varphi:V\to\{1,2,3\}$ is called a \textit{three-coloring} of $G$ if $\varphi(u)\neq\varphi(v)$ holds whenever $\{u,v\}\in E$.

\section{The reduction}

We proceed with the description of our reduction. We denote by $n,m$ the cardinalities of $V,E$, respectively, we set
\begin{equation}\label{deft}t=\frac{1}{4(m+n)^3},\end{equation}
and we denote by $V^i=\{v_1^i,\ldots,v_n^i\}$ a copy of the set $V=\{v_1,\ldots,v_n\}$. Our reduction from three-coloring to CSSP is the matrix $\mathcal{M}=\M(G)$ defined as follows. The rows of $\M$ have indexes in $V\cup\{1,2,3\}\cup\{\ep\}$, and the columns are indexed with $V^1\cup V^2\cup V^3\cup E$. We denote by $\M(\alpha,\beta)$ the entry of $\M$ at the intersection of a row with index $\alpha$ and a column with index $\beta$. For all distinct vertices $u,v\in V$, for all distinct numbers $i,j\in\{1,2,3\}$, and for all $e\in E$, we define the entries of $\M$ by

\noindent (1) $\M(u,v^i)=0$, $\M(u,u^i)=1$,

\noindent (2) $\M(u,e)=t^2$ if $u$ is adjacent to $e$ and $\M(u,e)=0$ otherwise,

\noindent (3) $\M(i,v^j)=0$, $\M(i,v^i)=t^3$, $\M(i,e)=t^5$,

\noindent (4) $\M(\ep,v^i)=0$, $\M(\ep,e)=t$.

For instance, the matrix $\M(K_3)$, which corresponds to the triangle graph $K_3$, looks like
$$
\left(\begin{array}{ccc|ccc|ccc|ccc}
1&0&0&1&0&0&1&0&0&t^2&0&t^2\\
0&1&0&0&1&0&0&1&0&0&t^2&t^2\\
0&0&1&0&0&1&0&0&1&t^2&t^2&0\\\hline
t^3&t^3&t^3&0&0&0&0&0&0&t^5&t^5&t^5\\
0&0&0&t^3&t^3&t^3&0&0&0&t^5&t^5&t^5\\
0&0&0&0&0&0&t^3&t^3&t^3&t^5&t^5&t^5\\\hline
0&0&0&0&0&0&0&0&0&t&t&t
\end{array}\right),
$$
where the first three rows have indexes in $V$, the next three rows have indexes $\{1,2,3\}$, and the last row has index $\ep$. Similarly, the columns of the above matrix are partitioned into triples corresponding to $V^1, V^2, V^3$, and $E$, respectively. Our goal in the rest of the paper is to show the following.

\begin{thr}\label{thr3cCSSP}
A graph $G$ has a three-coloring if and only if $$\delta_n\left(\M(G)\right)\leqslant \sqrt{mt^2+4nt^6+mt^{10}}.$$
\end{thr}

This theorem would complete the NP-hardness proof of Problem~\ref{prob11}.

\section{The proof. Part one}

In the rest of our paper, we denote by $\V$ an $n$-element subset of the set of column indexes of $\M$. The submatrix consisting of the columns of $\M$ with indexes in $\V$ is to be denoted by $S$. Also, we write $\M(\square,\beta)$ to denote the column of $\M$ with index $\beta$. As said above, the two conditions
\begin{equation}\label{eqnorm}\min\limits_{A} \|\M-SA\|\leqslant\sqrt{mt^2+4nt^6+mt^{10}},\end{equation}
\begin{equation}\label{eqnorm2}\|(I-SS^+)\M\|\leqslant\sqrt{mt^2+4nt^6+mt^{10}}\end{equation}
are equivalent, and Theorem~\ref{thr3cCSSP} states that one of them is satisfied by some $S$ if and only if $G$ is a three-colorable graph. In this section, we point out several necessary conditions for $S$ to satisfy~\eqref{eqnorm} and~\eqref{eqnorm2}.

\begin{lem}\label{lemcase2}
If the sets $\V$ and $E\cup \{v^1,v^2,v^3\}$ are disjoint for some $v\in V$, then the inequality~\eqref{eqnorm} is false.
\end{lem}

\begin{proof}
We see that the $v$th entry of any column of $S$ is zero, so the matrix $\M-SA$ has ones at $(v^1,v), (v^2,v), (v^3,v)$. We have $\|\M-SA\|\geqslant3$, which contradicts~\eqref{eqnorm}.
\end{proof}

In order to prove one more result of similar kind, we need to recall the definition of a \textit{strongly column diagonally dominant} matrix. Namely, this is a real square matrix $D$ any of whose diagonal elements is greater than the sum of the absolute values of all the elements in the same column. A classical result of matrix theory, sometimes called the \textit{Levy-Desplanques theorem}, states that any such matrix is non-singular. (An interested reader is referred to~\cite{Taus} for a historical survey on this theorem.)

\begin{lem}\label{lemcase1}
If there is $e\in\V\cap E$, then~\eqref{eqnorm} is false.
\end{lem}

\begin{proof}
Denoting by $\mathcal{S}$ the linear space spanned by the columns of $S$, we see that $\M(\square, e)$ belongs to $\mathcal{S}$. Taking into account~\eqref{deft} and~\eqref{eqnorm}, we conclude that, for all $v\in V$, the space $\mathcal{S}$ contains a column $C_v$ such that $\|C_v-\M(\square,v^1)\|<(n+1)^{-1}$. In particular, the $v$th element of $C_v$ is at least $n/(n+1)$, and the absolute value of the $u$th element of $C_v$ is less than $(n+1)^{-1}$ (where $u\in V\cup\{\varepsilon\}$ is different from $v$).

Since $\mathcal{S}$ is spanned by $n$ vectors, the matrix (we denote it by $C$) consisting of the $n+1$ columns $\{C_v\}$ and $\M(\square, e)$ must be of rank at most $n$. However, the submatrix of $C$ formed by taking the rows with indexes in $V\cup\{\varepsilon\}$ is strongly column diagonally dominant, so we get a contradiction.
\end{proof}

\begin{lem}\label{lemcase3}
The inequalities~\eqref{eqnorm},~\eqref{eqnorm2} imply that $\V=\{v^{\psi(v)},v\in V\}$, where $\psi$ is a mapping from $V$ to $\{1,2,3\}$.
\end{lem}

\begin{proof}
Follows immediately from Lemmas~\ref{lemcase2} and~\ref{lemcase1}.
\end{proof}

\section{The proof. Part two}

In this section, we consider the indexing sets $\V$ as in Lemma~\ref{lemcase3}. Namely, we will prove that the matrix $S$ satisfies~\eqref{eqnorm2} if and only if the corresponding mapping $\psi$ is a three-coloring of $G$. This would mean that, if $G$ does not admit a three-coloring, the inequality~\eqref{eqnorm2} is never satisfied. This would prove Theorem~\ref{thr3cCSSP} and show that the task of determining an optimal submatrix $S$ in a CSSP instance is at least as hard as finding a three-coloring of a given graph.

\begin{lem}\label{obs51}
If $G$ is three-colorable, then $\delta_n\left(\M(G)\right)\leqslant \sqrt{mt^2+4nt^6+mt^{10}}$.
\end{lem}

\begin{proof}
Let $\varphi: V\to \{1,2,3\}$ be a three-coloring of $G$. We take $$V^*=\{v^{\varphi(v)}, v\in V\}$$
to be the indexes of columns of $\M$ which form the matrix $S$ as in~\eqref{eqnorm}, and we proceed with the description of $A$. Its rows will have indexes in $V^*$, and the indexing of columns will correspond to that of $\M$. We set $A(\square,v^i)$ to be the column with the one at the position $v^{\varphi(v)}$ and zeros everywhere else; for $e=\{u,v\}$ in $E$, the column $A(\square,e)$ has $t^2$ at the positions $u^{\varphi(u)},v^{\varphi(v)}$ and zeros elsewhere. It remains to check that $\|\M-SA\|=\sqrt{mt^2+4nt^6+mt^{10}}$.
\end{proof}

\begin{lem}\label{lemcase34}
Assume $\V=\{v^{\psi(v)},v\in V\}$, where a mapping $\psi$ from $V$ to $\{1,2,3\}$ is not a three-coloring of $G$. Then~\eqref{eqnorm2} is not satisfied.
\end{lem}

\begin{proof}
We can write
$$
S=\left(\begin{array}{c|c|c}
I_{n_1}&O&O\\\hline
O&I_{n_2}&O\\\hline
O&O&I_{n_3}\\\hline
t^3 \mathbf{j}_{n_1}&O&O\\\hline
O&t^3 \mathbf{j}_{n_2}&O\\\hline
O&O&t^3 \mathbf{j}_{n_3}\\\hline
O&O&O
\end{array}\right),
$$
where the splitting of columns into blocks corresponds to the indexing sets $\psi^{-1}(1)$, $\psi^{-1}(2)$, $\psi^{-1}(3)$, $\{1\}$, $\{2\}$, $\{3\}$, $\{\ep\}$. The $O$'s stand for zero matrices of relevant sizes, and $\mathbf{j}_{k}$ is the $1\times k$ matrix of all ones. Since the columns of $S$ are linearly independent, we can write $S^+=(S^\top S)^{-1}S^\top$, and a computation shows that $I-SS^+$ equals
\begin{equation}\label{eqiss}
\left(\begin{array}{c|c|c|c|c|c|c}
u_1t^6J_{n_1}&O&O&-u_1t^3\mathbf{j}^\top_{n_1}&O&O&O\\\hline
O&u_2t^6J_{n_2}&O&O&-u_2t^3\mathbf{j}^\top_{n_2}&O&O\\\hline
O&O&u_3t^6J_{n_3}&O&O&-u_3t^3\mathbf{j}^\top_{n_3}&O\\\hline
-u_1t^3\mathbf{j}_{n_1}&O&O&u_1&0&0&0\\\hline
O&-u_2t^3\mathbf{j}_{n_2}&O&0&u_2&0&0\\\hline
O&O&-u_3 t^3\mathbf{j}_{n_3}&0&0&u_3&0\\\hline
O&O&O&0&0&0&1
\end{array}\right),
\end{equation}
where $u_i=({1+n_it^6})^{-1}$. Again, the splitting of the rows and columns into blocks corresponds to the partition $\psi^{-1}(1)$, $\psi^{-1}(2)$, $\psi^{-1}(3)$, $\{1\}$, $\{2\}$, $\{3\}$, $\{\ep\}$. 

By~\eqref{eqiss}, the $\ep$th row of $(I-SS^+)\M$ equals the $\ep$th row of $\M$, so it has $m$ entries equal to $t$.
Now let $\Delta$ be the matrix formed by the rows of $(I-SS^+)\M$ with indexes $1,2,3$. Using~\eqref{eqiss} once again, we note that, for $v\in V$ and $j\in\{1,2,3\}\setminus\{\psi(v)\}$, the $v^j$th column of $\Delta$ contains one zero, one element of the form $u_it^3$, and one element of the form $-u_it^3$. (Although it is not relevant for the proof, let us note that the $v^{\psi(v)}$th column of $\Delta$ consists of zeros.)

For $e=\{u,v\}\in E$, there are two possibilities. If $\psi(u)=\psi(v)$, then the $e$th column of $\Delta$ contains two zeros and one element of the form $u_it^5$. Otherwise, all the three elements in that column are of the form $\pm u_it^5$ (and since $\psi$ is not a three-coloring of $G$, this occurrence arises at least once). Therefore, the sum of squares of the entries of $(I-SS^+)\M$ is greater than or equal to
$$mt^2+4n\left(\frac{t^3}{1+nt^6}\right)^2+(m+2)\left(\frac{t^5}{1+nt^{6}}\right)^2,$$
which is greater than $mt^2+4nt^6+mt^{10}$ as a comparison with~\eqref{deft} shows.
\end{proof}

This completes the proof of the main results of our paper. Namely, we have shown that the column subset selection problem is NP-hard, and this problem becomes NP-complete in the Turing model of computation when restricted to rational matrices. Recall that we were working with the Frobenius norm of matrices, which seems to be the most common way to measure the cost function in the CSSP problem (see~\cite{BMD, Civ}). However, the same method will surely allow one to prove the same results for CSSP with respect to other norms like the $\ell_p$ operator norms for all $p>1$ or the spectral norm of matrices.

\end{document}